\documentclass[11pt,reqno]{amsart}
\usepackage{amsfonts}
\usepackage{amsmath,amsthm,amssymb,amsfonts,amscd}
\usepackage{mathrsfs}
\usepackage{lipsum}
\usepackage{bbding}
\usepackage{graphicx,latexsym}
\usepackage{hyperref}

\newcommand{\bea}{\begin{eqnarray}}
\newcommand{\eea}{\end{eqnarray}}
\newcommand{\bna}{\begin{eqnarray*}}
\newcommand{\ena}{\end{eqnarray*}}

\numberwithin{equation}{section}

\renewcommand{\thefootnote}{\fnsymbol{footnote}}
\theoremstyle{plain}
\newtheorem{theorem}{Theorem}
\newtheorem{lemma}{Lemma}

\newtheorem{proposition}{Proposition}

\theoremstyle{definition}

\newtheorem{remark}{Remark}

\newcommand\blfootnote[1]{%
  \begingroup
  \renewcommand\thefootnote{}\footnote{#1}%
  \addtocounter{footnote}{-1}%
  \endgroup
}

\begin{document}

\title{Bounds for $GL_3$ $L$-functions in depth aspect}

\author{Qingfeng Sun and Rui Zhao}

\begin{abstract}
Let $f$ be a Hecke-Maass cusp form for $SL_3(\mathbb{Z})$ and $\chi$
a primitive Dirichlet character of prime power conductor $\mathfrak{q}=p^{\kappa}$ with $p$ prime
and $\kappa\geq 10$. We prove a subconvexity bound
$$
L\left(\frac{1}{2},\pi\otimes \chi\right)\ll_{p,\pi,\varepsilon}
\mathfrak{q}^{3/4-3/40+\varepsilon}
$$
for any $\varepsilon>0$, where the dependence of the implied constant on $p$ is explicit and polynomial.
We obtain this result by applying the circle method of Kloosterman's version,
summation formulas of Poisson and Voronoi's type and a conductor lowering mechanism
introduced by Munshi \cite{Munshi31}. The main new technical estimates are the essentially
square root bounds for some twisted multi-dimensional character sums, which are proved by
an elementary method.
\end{abstract}


\maketitle
\tableofcontents

\section{Introduction}

Let $L(s,f)$ be an $L$-function with the analytic conductor $\mathfrak{q}(s,f)$.
By the functional equation and the Phragmen-Lindel\"{o}f convexity principle,
we have the convexity bound $L(s,f)\ll \mathfrak{q}(s,f)^{1/4+\varepsilon}$.
It is an fascinating problem to break the convexity barrier. In the $t$-aspect,
one has the classical result for the Riemann zeta function
$\zeta(1/2+it)\ll_{\varepsilon} (1+|t|)^{1/6+\varepsilon}$ due
to Weyl \cite{W}. For $L$-functions on $GL_2$, results of the same strength
\bea
L\left(\frac{1}{2}+it,f\right)\ll_{f,\varepsilon} (1+|t|)^{1/3+\varepsilon}
\eea
were proved by Good \cite{Good}, Jutila \cite{Jutila} and Meurman \cite{Meur},
where $f$ is a fixed holomorphic cusp form or a Maass
cusp form. For $GL_3$ $L$-functions, Munshi \cite{Munshi32} proved that
\bea
L\left(\frac{1}{2}+it,\pi\right)\ll_{\pi,\varepsilon}
(1+|t|)^{3/4-1/16+\varepsilon},
\eea
where $\pi$ is a fixed $GL_3$ Hecke-Maass cusp form (this bound was first proved by
Li \cite{Li} for $\pi$ self-dual).
On the other hand, in the conductor aspect,
we have the Burgess' bounds $L(1/2,\chi)\ll_{\varepsilon} q^{3/16+\varepsilon}$
and $L\left(1/2,f\otimes\chi\right)\ll_{f,\varepsilon}q^{3/8+\varepsilon}$
for a primitive character $\chi$ of conductor $q$, where
$f$ is a fixed $GL_2$ cusp form. Interestingly, for $\chi$ quadratic,
Conrey and Iwaniec \cite{CI} proved the exponent $1/3$, i.e., the quantitative
analogue of (1.1). Recently, by developing
a general result on $p$-adic analytic phase and a $p$-adic version of van der Corput's
method for exponential sums, Blomer and Mili\'{c}evi\'{c} \cite{BM} also
proved the same exponent if the conductor of $\chi$ is a prime power $\mathfrak{q}=p^{\kappa}$
\bea
L\left(1/2+it,f\otimes\chi\right)\ll_{p,t,f,\varepsilon}\mathfrak{q}^{1/3+\varepsilon},
\eea
where the implied constant on $p$ and $t$ is explicit and polynomial (Munshi and Singh proved
the same result using the approach in \cite{Munshi31}). Also see
\cite{Mili} and \cite{Le} for other interesting subconvexity results in the depth aspect.

Let $\pi$ be a Hecke-Maass cusp form for $SL_3(\mathbb{Z})$ and $\chi$
a primitive Dirichlet character modulo $\mathfrak{q}$.
Then the convexity bound for $L\left(1/2,\pi\otimes \chi\right)$ is $\mathfrak{q}^{3/4+\varepsilon}$.
For $\mathfrak{q}$ prime, the subconvexity results for $L\left(1/2,\pi\otimes \chi\right)$
have recently been established
in the work \cite{B}, \cite{HN}
and \cite{Munshi4}-\cite{Munshi5}. Munshi \cite{Munshi31} showed
a subconvexity bound for $\mathfrak{q}$ square-free.
In this paper, following Munshi \cite{Munshi31}, we want to prove a subconvexity bound
for $L\left(1/2,\pi\otimes \chi\right)$
in the depth aspect. Our main result is the following.

\begin{theorem}
Let $\pi$ be a Hecke-Maass cusp form for $SL_3(\mathbb{Z})$ and $\chi$
a primitive Dirichlet character of prime power conductor $\mathfrak{q}=p^{\kappa}$ with $\kappa\geq 3$. We have
$$
L\left(\frac{1}{2},\pi\otimes \chi\right)\ll_{\pi,\varepsilon}
p^{3/4}
\mathfrak{q}^{3/4-3/40+\varepsilon}
$$
for any $\varepsilon>0$.
\end{theorem}

\begin{remark}
Our result in Theorem 1 can be compared with the $t$-aspect subconvexity in (1.2)
as explained in \cite{Munshi2}. It is worth noting that for $\pi$ the symmetric-square lifts of $GL_2$ cusp forms,
Munshi \cite{Munshi2} proved the better result $\mathfrak{q}^{3/4-1/12+\varepsilon}$
by the moment method.
\end{remark}

\begin{remark}
We are not trying to get the best exponent in $p$.
With the present exponent $3/4$, the bound in Theorem 1 breaks the convexity for $\kappa>10$.
\end{remark}

\medskip
\noindent
{\bf Notation.}
Throughout the paper, the letters $q$, $m$ and $n$, with or without subscript,
denote integers. The letter $\varepsilon$ is an arbitrarily small
positive constant, not necessarily the same at different occurrences. The symbol
$\ll_{a,b,c}$ denotes that the implied constant depends at most on $a$, $b$ and $c$.
 Finally, fractional numbers
such as $\frac{ab}{cd}$ will be written as $ab/cd$ and $a/b+c$ or $c+a/b$ means
$\frac{a}{b}+c$.

\section{Sketch of the proof}

By the functional equation we have $L\left(\frac{1}{2},\pi\otimes \chi\right)\ll N^{-1/2}\mathscr{S}(N)$, where
\bna
\mathscr{S}(N)=\sum_{n\sim N}A_{\pi}(1,n)\chi(n),
\ena
with $N\sim \mathfrak{q}^{3/2}$. Applying the conductor lowering mechanism
introduced by Munshi \cite{Munshi31}, we have
\bna
\mathscr{S}(N)=\mathop{\sum\sum}_{n,m\sim N\atop n\equiv
m ({\rm mod}p^{\lambda})}A_{\pi}(1,n)\chi(n)\delta\left(\frac{n-m}{p^\lambda}\right)
\ena
where $\delta: \mathbb{Z}\rightarrow \{0,1\}$ with
$\delta(0)=1$ and $\delta(n)=0$ for $n\neq 0$, and $\lambda\geq 2$ is an integer to be chosen later.
Using Kloosterman's circle method and removing the congruence
$n\equiv
m ({\rm mod}p^{\lambda})$ by exponential sums we get
\bna
\mathscr{S}(N)\approx \frac{1}{p^{\lambda}}\sum_{q\sim Q}\frac{1}{q}
\sum_{Q<a\leq q+Q \atop (a,q)=1} \frac{1}{a}
   \sum_{b(\text{mod} \,p^{\lambda})}\mathop{\sum\sum}_{n,m\sim N}
   A_{\pi}(1,n)\chi(m)e\left(\frac{(\overline{a}+bq)(n-m)}{qp^\lambda}\right).
\ena
Trivially we have $\mathscr{S}(N)\ll N^2$.

For simplicity, we assume $(q,p)=1$ and $(\overline{a}+bq,p)=1$.
Recall $\chi$ is of modulus $\mathfrak{q}=p^{\kappa}$.
Then the conductor of the $m$-sum has the size $qp^{\kappa}$. Applying Poisson summation
to the $m$-sum we get that the dual sum is of size $qp^{\kappa}/N$.
The conductor for the $n$-sum has the size $qp^\lambda$ and the dual sum after
$GL_3$ Voronoi summation formula
is essentially supported on summation of size $q^3p^{3\lambda}/N$.
Assuming
square-root cancellation for the character sum,
we find that we have saved
\bna
\frac{N}{(qp^{\kappa})^{1/2}}\times \frac{N}{(qp^{\lambda})^{3/2}}\times (qp^{\lambda})^{1/2}
\sim \mathfrak{q}^{11/8}p^{-\lambda/4}.
\ena

Now we arrive at an expression
of the form
\bna
\sum_{1\leq n\ll Q^3p^{3\lambda}/N}A_{\pi}(n,1)
\sum_{q\sim Q}\chi(q)
\sum_{|m|\ll Qp^{\kappa}/N}\quad
\sum_{b({\rm mod} p^{\lambda})}\overline{\chi}(m-bp^{\kappa-\lambda})
S(\overline{b},n;qp^{\lambda}).
\ena
Next we apply Cauchy-Schwartz inequality to get rid of the Fourier coefficients.
Then we need to deal with
\bna
\sum_{1\leq n\ll Q^3p^{3\lambda}/N}\left|
\sum_{q\sim Q}\chi(q)
\sum_{|m|\ll Qp^{\kappa}/N}\quad
\sum_{b({\rm mod} p^{\lambda})}\overline{\chi}(m-bp^{\kappa-\lambda})
S(\overline{b},n;qp^{\lambda})\right|^2.
\ena
Opening the square and applying Poisson summation to the
sum over $n$, we are able to save $Q^2p^{\kappa}/N\sim p^{\kappa-\lambda}$ from the diagonal term
and $$\frac{Q^3p^{3\lambda}/N}{\sqrt{Q^2p^{\lambda}}}\sim p^{3\lambda/2}$$ from the off-diagonal term. So the optimal choice for $\lambda$
is given by $\lambda=2\kappa/5$. In total, we have saved
\[
\mathfrak{q}^{11/8}p^{-\lambda/4}\times p^{3\lambda/4}\sim \mathfrak{q}^{3/2+3/40}.
\]
It follows that
\bna
L\left(\frac{1}{2},\pi\otimes \chi\right)\ll N^{-1/2}\mathscr{S}(N)\ll N^{3/2}\mathfrak{q}^{-3/2-3/40}\sim \mathfrak{q}^{3/4-3/40}.
\ena

\section{Proof of Theorem 1}

By the approximate functional equation we have
\bea
L\left(\frac{1}{2},\pi\otimes \chi\right)\ll_{\pi,\varepsilon} \mathfrak{q}^{\varepsilon}
\sup_{N\leq \mathfrak{q}^{3/2+\varepsilon}}\frac{|\mathscr{S}(N)|}{\sqrt{N}},
\eea
where
\bna
\mathscr{S}(N)=\sum_{n}A_{\pi}(1,n)\chi(n)
V\left(\frac{n}{N}\right)
\ena
for some smooth function $V$ supported in $[1,2]$ and satisfying $V^{(j)}(y)\ll_j 1$.
Note that by Cauchy's inequality and the Rankin-Selberg estimate (see \cite{Mol})
\bea
\mathop{\sum\sum}_{n_1^2n_2\leq Y} \left|A_{\pi}(n_1,n_2)\right|^2\ll_{\pi,\varepsilon}
Y^{1+\varepsilon},
\eea
we have the trivial bound $\mathscr{S}(N)\ll_{\pi,\varepsilon}N$. Thus Theorem
1 is true for
$N\leq \mathfrak{q}^{27/20}$.
In the following, we will estimate $\mathscr{S}(N)$ in the range
\bea
\mathfrak{q}^{27/20}<N\leq \mathfrak{q}^{3/2+\varepsilon}.
\eea

\begin{proposition}
Assume $\lambda\leq 2\kappa/3$ and (3.3). Then we have
\bna
\mathscr{S}(N)\ll N^{1/2+\varepsilon}(p^{3\kappa/8+3\lambda/4}
+p^{7\kappa/8-\lambda/2+3/4}).
\ena
\end{proposition}
Take $\lambda=\lfloor2\kappa/5\rfloor+1$, where $\lfloor x\rfloor$ denotes
the largest integer which does not exceed $x$.
By (3.3) and Proposition 1, we have
\bna
\mathscr{S}(N)\ll p^{3/4}N^{1/2+\varepsilon}\mathfrak{q}^{3/4-3/40}
\ena
Then Theorem 1 follows from above bound and (3.1).
In the following we prove Proposition 1.

\subsection{The circle method}
Define $\delta: \mathbb{Z}\rightarrow \{0,1\}$ with
$\delta(0)=1$ and $\delta(n)=0$ for $n\neq 0$.
By Kloosterman's version
of the circle method, for any $n\in \mathbb{Z}$ and $Q\in \mathbb{R}^+$, we have
\bea
\delta(n)=2\mathrm{Re}\int_0^1 \sum_{1\leq q\leq Q}\sum_{Q<a\leq q+Q \atop (a,q)=1}
\frac{1}{aq}e\left(\frac{n\overline{a}}{q}-\frac{n\zeta}{aq}\right)\mathrm{d}\zeta,
\eea
where throughout the paper $e(z)=e^{2\pi iz}$ and $\overline{a}(\mathrm{mod}\,q)$ denotes
the multiplicative inverse of $a$ modulo $q$. Define $\mathbf{1}_\mathscr{F}=1$ if $\mathscr{F}$ is true, and is 0 otherwise.
Following Munshi \cite{Munshi31} we write $\delta(n)$ as
$\delta(n/p^\lambda)\mathbf{1}_{p^\lambda|n}$ ($2\leq \lambda<\kappa, \lambda\in \mathbb{N}$
is a parameter to be determined later)
to lower the conductor and obtain
\bna
\mathscr{S}(N)=\sum_{n}
A_{\pi}(1,n)V\left(\frac{n}{N}\right)
\sum_{p^\lambda|n-m}
\chi(m)U\left(\frac{m}{N}\right)
\delta\left(\frac{n-m}{p^\lambda}\right),
\ena
where $U$ is a smooth function supported in $[1/2,5/2]$, $U(y)=1$
for $y\in [1,2]$ and $U^{(j)}(y)\ll_j 1$.
Applying (3.4) and choosing
\bna
Q=\sqrt{N/p^\lambda}
\ena
we get
\bna
\mathscr{S}(N)=\mathscr{S}^+(N)+\mathscr{S}^-(N),
\ena
where
\bna
\mathscr{S}^{\pm}(N)&=&\int_0^1\sum_{1\leq q\leq Q}\sum_{Q<a\leq q+Q\atop (a,q)=1}
\frac{1}{aq}
\sum_{n}
A_{\pi}(1,n)V\left(\frac{n}{N}\right)
\\&&\sum_{p^\lambda|n-m}
\chi(m)U\left(\frac{m}{N}\right)
e\left(\pm\frac{(n-m)\overline{a}}{qp^\lambda}\mp\frac{(n-m)\zeta}{aqp^\lambda}\right)\mathrm{d}\zeta.
\ena

We will only estimate $\mathscr{S}^+(N)$ (the same analysis
holds for $\mathscr{S}^-(N)$) and write $\mathscr{S}^+(N)$ as $\mathscr{S}(N)$. Removing the condition $p^\lambda|n-m$ using exponential sums to separate the variables $m$ and $n$
we get
\bea
\mathscr{S}(N)=\int_{0}^{1}\sum_{1\leq q\leq Q} \sum_{Q<a\leq q+Q \atop (a,q)=1} \frac{1}{aqp^{\lambda}}
   \sum_{b(\text{mod} \,p^{\lambda})} \mathscr{A}\times \mathscr{B} \hspace{3pt} \mathrm{d}\zeta,
\eea
where
$$
\mathscr{A}=\sum_{m}\chi(m)
   e\left(-\frac{(\overline{a}+bq)m}{qp^{\lambda}}\right)U\left(\frac{m}{N}\right)
   e\left(\frac{m\zeta}{aqp^{\lambda}}\right)
$$
and
$$
\mathscr{B}=\sum_{n}A_{\pi}(1,n) e\left(\frac{(\overline{a}+bq)n}{qp^{\lambda}}\right)
  V\left(\frac{n}{N}\right)e\left(-\frac{n\zeta}{aqp^{\lambda}}\right).
$$

\subsection{Summation formulas and Cauchy-Schwartz}

Next we transform $\mathscr{A}$ and $\mathscr{B}$ by Poisson summation formula and $GL_3$
Voronoi formula, respectively, and obtain the following results.
\begin{lemma}
Let $q=p^sq'$, $(q',p)=1$ and $s\geq 0$.
Then we have
$$
\mathscr{A}=\frac{N\chi(q')\tau_{\chi} }{p^{\kappa}} \sum_{|m|\leq N^{\varepsilon}Qp^{\kappa}/N
\atop m\equiv\overline{a}p^{\kappa-\lambda}(\mathrm{mod}\, q)}
\overline{\chi}\left(\frac{m-(\overline{a}+bq)p^{\kappa-\lambda}}{p^{s}}\right)
\mathfrak{I}(m,a,q,\zeta)+O(\mathfrak{q}^{-A})
$$
for any $A>0$, where the integral $\mathfrak{I}(m,a,q,\zeta)$ is defined in (4.2).
\end{lemma}

\begin{lemma}
Let $a^*=(\overline{a}+bq)/(\overline{a}+bq,qp^{\lambda})$ and
$q^*=qp^{\lambda}/(\overline{a}+bq,qp^{\lambda})$.
Then we have
\bna
\mathscr{B}&=&\frac{N^{1/2}}{q^{*1/2}}\sum_{\pm}\sum_{n_{1}|q^*}\,
\mathop{\sum\sum}_{n_1^2n_2\leq N^{\varepsilon}q^{*3}Q^3/q^3N}\frac{A_{\pi}(n_{2},n_{1})}{\sqrt{n_2}}
  S\left(\overline{a^*},\pm n_{2};\frac{q^*}{n_{1}}\right)\nonumber\\&&\qquad\qquad\qquad\qquad
  \qquad\qquad\qquad
\times\mathfrak{J}^{\pm}\left(\frac{n_{1}^{2}n_{2}}{q^{*3}},a,q,\zeta\right)+O(\mathfrak{q}^{-A})
\ena
for any $A>0$, where $\mathfrak{J}^{\pm}\left(y,a,q,\zeta\right)$ is defined in (5.2) and satisfies
\bna
\mathfrak{J}^{\pm}
\left(y,a,q,\zeta\right)\ll N^{\varepsilon}\sqrt{\frac{Q}{q}}.
\ena

\end{lemma}
The details of the proof of Lemmas 1 and 2 are in Sections 4 and 5.
Note that for  $s\geq 1$, we have $(\overline{a}+bq, qp^{\lambda})=1$, $a^*=\overline{a}+bq$ and
$q^*=qp^{\lambda}$. For $s=0$, we have
$(\overline{a}+bq, qp^{\lambda})=p^{r}$, $0\leq r\leq\lambda$, $a^*=(\overline{a}+bq)/p^r$ and
$q^*=qp^{\lambda-r}$. Since $\overline{a}+bq\equiv 0 (\text{{\rm mod }} p^r)$, we have
$b\equiv -\overline{aq} (\text{{\rm mod }} p^r)$. Denote $\varpi_q^r:=(1-q\overline{q})/p^r\in \mathbb{Z}$. We write
$\overline{a}+bq=(\overline{a}\varpi_q^r+cq)p^r$ with $c(\text{{\rm mod }} p^{\lambda-r})$.
Plugging Lemmas 1 and 2 into (3.5) and reducing
the $n_1, n_2$ sums into dyadic intervals, we have
\bea
\mathscr{S}(N)&\ll& \sum_{\pm}\sum_{r=0}^{\lambda}\sum_{L_1\ll N^{\varepsilon}p^{3\lambda-3r}Q^3/N\atop L_1 \, \mathrm{dyadic}}
|\mathscr{S}_1^{\pm}(N,L_1,r)|\nonumber\\
&&+
\sum_{\pm}\sum_{s=1}^{\log Q/\log p}
\sum_{L_2\ll N^{\varepsilon}p^{3\lambda}Q^3/N\atop
L_2 \, \mathrm{dyadic}}|\mathscr{S}_2^{\pm}(N,L_2,s)|+\mathfrak{q}^{-2018},
\eea
where
\bea
\mathscr{S}_1^{\pm}(N,L_1,r)&=&\frac{N^{3/2}}{p^{(\kappa+3\lambda-r)/2}}
   \mathop{\sum\sum}_{L_1/2<n_1^2n_2\leq L_1} \frac{A_{\pi}(n_2,n_1)}{\sqrt{n_2}}
 \sum_{1\leq q\leq Q, (q,p)=1\atop n_{1}|qp^{\lambda-r}} \frac{\chi(q)}{q^{3/2}}
   \nonumber\\&&
 \times
   \sum_{1\leq|m|\leq N^{\varepsilon}Qp^{\kappa}/N
   \atop (m,q)=1}
 \sum_{Q<a\leq q+Q \atop a\equiv\overline{m}p^{\kappa-\lambda}(\text{{\rm mod }} q)} \frac{1}{a}
   \mathfrak{K}^{\pm}\left(m,\frac{n_1^2n_2}{q^3p^{3\lambda-3r}},a,q\right)\nonumber\\
&&\times\sum_{c(\text{{\rm mod }} p^{\lambda-r})}
\overline{\chi}\left(m-cp^{\kappa-\lambda+r}\right)
   S\left(\overline{c},\pm n_2;\frac{qp^{\lambda-r}}{n_1}\right)\nonumber\\
\eea
and
\bna
\mathscr{S}_2^{\pm}(N,L_2,s)&=&\frac{N^{3/2}}{p^{(\kappa+3\lambda+3s)/2}}
   \mathop{\sum\sum}_{L_2/2<n_1^2n_2\leq L_2} \frac{A_{\pi}(n_2,n_1)}{\sqrt{n_2}}
   \sum_{1\leq q\leq Q/p^s, (q,p)=1 \atop n_1|qp^{\lambda+s}} \frac{\chi(q)}{q^{3/2}}
   \nonumber\\
&& \sum_{1\leq |m|\leq N^{\varepsilon}Qp^{\kappa}/N \atop
m\equiv0(\text{{\rm mod }} p^s)}
\sum_{Q<a\leq qp^s+Q, (a,p)=1\atop
a\equiv\overline{m}p^{\kappa-\lambda}(\text{{\rm mod }} q)} \frac{1}{a}
\mathfrak{K}^{\pm}\left(m,\frac{n_1^2n_2}{q^3p^{3\lambda+3s}},a,qp^s\right)\nonumber\\
&&\times\sum_{b(\text{{\rm mod }} p^{\lambda})}
\overline{\chi}\left(\frac{m-(\overline{a}+bp^s)p^{\kappa-\lambda}}{p^s}\right)
   S\left(\overline{\overline{a}+bp^s},\pm n_{2};\frac{qp^{\lambda+s}}{n_1}\right)
\ena
with
\bea
\mathfrak{K}^{\pm}(y_1,y_2,a,q)=\int_0^1\mathfrak{I}(y_1,a,q,\zeta)
\mathfrak{J}^{\pm}\left(y_2,a,q,\zeta\right)\mathrm{d}\zeta\ll N^{\varepsilon}\sqrt{\frac{q}{Q}}.
\eea
Here we have changed variables $\overline{a}\varpi_q^r+cq\rightarrow c$ and $bq\rightarrow b$.

\begin{remark}
If $m=0$, then the conditions $p^{\kappa-\lambda}\equiv 0 (\text{{\rm mod }} qp^s)$
and $((\overline{a}+bqp^s)p^{\kappa-\lambda-s},p)=1$
imply that $q=1$ and $s=\kappa-\lambda$. Thus we have
$p^{\kappa-\lambda}\leq Q=\sqrt{N/p^{\lambda}}$ which implies $N> p^{(3/2+\varepsilon)\kappa}$
which contradicts to the assumption (3.3). Therefore, we have $m\neq 0$.
\end{remark}

Applying Cauchy-Schwartz inequality to $n_1,n_2$-sums in (3.7) and using the Rankin-Selberg bound (3.2), we get
\bea
\mathscr{S}_1^{\pm}(N,L_1,r)&\ll&\frac{N^{3/2}L_1^{1/2}}{p^{(\kappa+3\lambda-r)/2}}\mathscr{H}_1^{\pm}(N,L_1,r)^{1/2}
\eea
where
\bea
   \mathscr{H}_1^{\pm}(N,L_1,r)&=&\sum_{(n_1',p)=1}\sum_{n_1''|p^{\lambda-r}}
   \sum_{n_2}\frac{1}{n_2}W\left(\frac{n_1'^2n_1''^2n_2}{L_1}\right)\left|
 \sum_{1\leq q\leq Q, (q,p)=1\atop n_1'|q} \frac{\chi(q)}{q^{3/2}}
  \sum_{1\leq|m|\leq N^{\varepsilon}Qp^{\kappa}/N
   \atop (m,q)=1}
 \right.\nonumber\\&&
 \left.
 \sum_{Q<a\leq q+Q \atop a\equiv\overline{m}p^{\kappa-\lambda}(\text{{\rm mod }} q)} \frac{1}{a}
\mathfrak{C}_r(m,n_1',n_1'',\pm n_2,a,q)
 \mathfrak{K}^{\pm}\left(m,\frac{n_1'^2n_1''^2n_2}{q^3p^{3\lambda-3r}},a,q\right)
\right|^2
\eea
with $W(y)$ a smooth positive function, $W(y)=1$ if $y\in [1/2,1]$, and
\bea
\mathfrak{C}_r(m,n_1',n_1'',n_2,a,q)=
S\left(a\overline{\varpi_q^r\widehat{p_r}},n_2\overline{\widehat{p_r}};\widehat{q}\right)
\sum_{c(\text{{\rm mod }} p^{\lambda-r})} \overline{\chi}\left(m-cp^{\kappa-\lambda+r}\right)
S\left(\overline{c}\overline{\widehat{q}},n_2\overline{\widehat{q}};\widehat{p_r}\right).
\eea
Here
$\widehat{q}=q/n_{1}'$ and $\widehat{p_r}=p^{\lambda-r}/n_{1}''$.
Similarly,
\bea
\mathscr{S}_2^{\pm}(N,L_2,s)\ll\frac{N^{3/2}L_2^{1/2}}{p^{(\kappa+3\lambda+3s)/2}}\mathscr{H}_2^{\pm}(N,L_2,s)^{1/2}
\eea
where
\bna
   \mathscr{H}_2^{\pm}(N,L_2,s)&=&\sum_{(n_1',p)=1}\sum_{n_1''|p^{\lambda+s}}\sum_{n_2}\frac{1}{n_2}
   W\left(\frac{n_1'^2n_1''^2n_2}{L_2}\right)\left|
   \sum_{1\leq q\leq Q/p^s, (q,p)=1 \atop n_1'|q} \frac{\chi(q)}{q^{3/2}}
   \sum_{1\leq |m|\leq N^{\varepsilon}Qp^{\kappa}/N \atop
m\equiv0(\text{{\rm mod }} p^s)}
\right.\nonumber\\
&&\left.\sum_{Q<a\leq qp^s+Q, (a,p)=1\atop
a\equiv\overline{m}p^{\kappa-\lambda}(\text{{\rm mod }} q)} \frac{1}{a}
\mathfrak{B}_s(m,n_1',n_1'',\pm n_2,a,q)
\mathfrak{K}^{\pm}\left(m,\frac{n_1'^2n_1''^2n_2}{q^3p^{3\lambda+3s}},a,qp^s\right)
\right|^2
\ena
with
\bea
\mathfrak{B}_s(m,n_1',n_1'',n_2,a,q)&=&
S\left(a\overline{p^{\lambda+s}/n_{1}''},n_2\overline{p^{\lambda+s}/n_{1}''};\frac{q}{n_{1}'}\right)
\sum_{b(\text{{\rm mod }} p^{\lambda})} \overline{\chi}\left(\frac{m-(\overline{a}+bp^s)p^{\kappa-\lambda}}{p^s}\right)\nonumber\\
   &&\qquad\qquad\qquad
   \times
   S\left(\overline{\overline{a}+bp^s}\,\overline{q/n_{1}'},n_2\overline{q/n_{1}'};\frac{p^{\lambda+s}}{n_{1}''}\right).
\eea

\subsection{Poisson summation}

Opening the square in (3.10) and switching the order of summations, we get

\bna
  \mathscr{H}_1^{\pm}(N,L_1,r)&=&\sum_{n_1'}\sum_{n_1''|p^{\lambda-r}}
  \sum_{1\leq q_1\leq Q,(q_{1},p)=1 \atop n_1'|q_1} \frac{\chi(q_1)}{q_1^{3/2}}
 \sum_{1\leq q_2\leq Q,(q_{2},p)=1\atop n_1'|q_2} \overline{\frac{\chi(q_2)}{q_2^{3/2}}}\\&&
     \sum_{1\leq|m_1|\leq N^{\varepsilon}Qp^{\kappa}/N \atop (m_1,q_1)=1}
    \sum_{1\leq|m_2|\leq N^{\varepsilon}Qp^{\kappa}/N\atop (m_2,q_2)=1 }
    \sum_{Q<a_1\leq q_1+Q
      \atop a_1\equiv\overline{m_1}p^{\kappa-\lambda}(\text{{\rm mod }} q_1)} \frac{1}{a_1}
   \sum_{Q<a_2\leq q_2+Q
   \atop a_2\equiv\overline{m_2}p^{\kappa-\lambda}(\text{{\rm mod }} q_2)} \frac{1}{a_2}\times
\mathbf{T},
\ena
where temporarily,
\bna
\mathbf{T}&=& \sum_{n_2}\frac{1}{n_2}W\left(\frac{n_1'^2n_1''^2n_2}{L_1}\right)
\mathfrak{C}_r(m_1,n_1',n_1'',\pm n_2,a_1,q_1)\overline{\mathfrak{C}_r(m_2,n_1',n_1'',\pm n_2,a_2,q_2)}\\&&\qquad\qquad\qquad\qquad\quad
  \mathfrak{K}^{\pm}\left(m_1,\frac{n_1'^2n_1''^2n_2}{q_1^3p^{3\lambda-3r}},a_1,q_1\right)
\overline{  \mathfrak{K}^{\pm}\left(m_2,\frac{n_1'^2n_1''^2n_2}{q_2^3p^{3\lambda-3r}},a_2,q_2\right)}.
\ena
Applying Poisson summation with modulus $\widehat{q_1}\widehat{q_2}\widehat{p_r}$, we arrive at
\bna
\mathbf{T}
=\frac{1}{\widehat{q_{1}}\widehat{q_{2}}\widehat{p_r}}  \sum_{n_{2}\in\mathbb{Z}} \mathfrak{C}^*\times \mathfrak{K}^*,
\ena
where
\bea
\mathfrak{C}^*=\sum_{\beta(\text{{\rm mod }}\widehat{q_1}\widehat{q_2}\widehat{p_r})}
\mathfrak{C}_r\left(m_{1},n_{1}',n_{1}'',\beta,a_{1},q_{1}\right)
      \overline{\mathfrak{C}_r\left(m_{2},n_1',n_1'',\beta,a_2,q_{2}\right)}
      e\left(\frac{\pm n_2\beta}{\widehat{q_1}\widehat{q_2}\widehat{p_r}}\right)
\eea
and
\bna
\mathfrak{K}^*=\int_{\mathbb{R}} W(y)  \mathfrak{K}^{\pm}\left(m_1,\frac{L_1y}{q_1^3p^{3\lambda-3r}},a_1,q_1\right)
\overline{\mathfrak{K}^{\pm}\left(m_2,\frac{L_1y}{q_2^3p^{3\lambda-3r}},a_2,q_2\right)}
e\left(-\frac{n_2L_1y}{q_1q_2p^{\lambda-r}n_1''}\right)\frac{\mathrm{d}y}{y}.
\ena
The integral $\mathfrak{K}^*$ gives arbitrary power saving in $\mathfrak{q}$ if
$|n_2|\geq N^{\varepsilon}Q^2p^{\lambda-r}n_1''/L_1$ for any $\varepsilon>0$.
For small values of $n_2$, by (3.8), we have
$$
\mathfrak{K}^*\ll N^{\varepsilon}\frac{\sqrt{q_{1}q_{2}}}{Q}.
$$
Therefore, at the cost of a negligible error,
\bna
\mathbf{T}
\ll N^{\varepsilon}\frac{\sqrt{q_{1}q_{2}}}{Q}\frac{1}{\widehat{q_{1}}\widehat{q_{2}}\widehat{p_r}}
\sum_{|n_2|\leq N^{\varepsilon}Q^2p^{\lambda-r}n_1''/L_1} |\mathfrak{C}^*|
\ena
and
\bea
  \mathscr{H}_1^{\pm}(N,L_1,r)&\ll&N^{\varepsilon}\frac{1}{Q^3}\sum_{n_1'}\sum_{n_1''|p^{\lambda-r}}
  \sum_{1\leq q_1\leq Q \atop n_1'|q_1} \frac{1}{q_1}
 \sum_{1\leq q_2\leq Q\atop n_1'|q_2} \frac{1}{q_2}\nonumber\\&&
  \sum_{1\leq|m_1|\leq N^{\varepsilon}Qp^{\kappa}/N \atop (m_1,q_1)=1}
    \sum_{1\leq|m_2|\leq N^{\varepsilon}Qp^{\kappa}/N\atop (m_2,q_2)=1 }
 \frac{1}{\widehat{q_{1}}\widehat{q_{2}}\widehat{p_r}}
\sum_{|n_2|\leq N^{\varepsilon}Q^2p^{\lambda-r}n_1''/L_1} |\mathfrak{C}^*|.
\eea
Similarly,
\bea
  \mathscr{H}_2^{\pm}(N,L_2,s)&\ll&N^{\varepsilon}\frac{p^{3s}}{Q^3}\sum_{n_1'}\sum_{n_1''|p^{\lambda+s}}
  \sum_{1\leq q_1\leq Q/p^s,(q_{1},p)=1 \atop n_1'|q_1} \frac{1}{q_1}
 \sum_{1\leq q_2\leq Q/p^s,(q_{2},p)=1\atop n_1'|q_2} \frac{1}{q_2}\nonumber\\&&
          \sum_{1\leq |m_1|\leq N^{\varepsilon}Qp^{\kappa}/N\atop
         m_1\equiv0(\text{{\rm mod }} p^s)}
    \sum_{1\leq |m_2|\leq N^{\varepsilon}Qp^{\kappa}/N\atop
    m_2\equiv0(\text{{\rm mod }} p^s)}
 \frac{1}{\widehat{q_{1}}\widehat{q_{2}}\widehat{\rho_s}}
\sum_{|n_2|\leq N^{\varepsilon}Q^2p^{\lambda-s}n_1''/L_2} |\mathfrak{B}^*|,
\eea
where $\widehat{q}=q/n_1'$, $\widehat{\rho_s}=p^{\lambda+s}/n_1''$ and
\bea
\mathfrak{B}^*=\sum_{\beta(\text{{\rm mod }}\widehat{q_1}\widehat{q_2}\widehat{\rho_s})}
 \mathfrak{B}_s\left(m_{1},n_{1}',n_{1}'',\beta,a_{1},q_{1}\right)
      \overline{\mathfrak{B}_s\left(m_{2},n_1',n_1'',\beta,a_2,q_{2}\right)}
      e\left(\frac{\pm n_2\beta}{\widehat{q_1}\widehat{q_2}\widehat{\rho_s}}\right).
\eea

\begin{lemma}
Assume $\lambda\leq 2\kappa/3$. Let $p^k\| n_2$ with $k\geq 0$.

(1) We have
\bna
\mathfrak{C}^*\ll \widehat{q_1}\widehat{q_2}(\widehat{q_1},\widehat{q_2},n_2)
\widehat{p_r}^2p^{2(\lambda-r)}.
\ena
Moreover, for $n_2=0$, the character sums vanish unless $q_1=q_2$ in which case
\bna
\mathfrak{C}^*\ll \widehat{q_1}^2(\widehat{q_1},m_1-m_2)\widehat{p_r}p^{2(\lambda-r)}.
\ena

(2) For $n_1''=p^{\lambda-r}$ or $n_1''=p^{\lambda-r-1}$ with $\lambda-r\geq 2$, we have
\bna
\mathfrak{C}^*=0.
\ena

(3) For $p^{\lambda-r}/n_1''\geq p^2$, we have
$\mathfrak{C}^*$ vanishes unless $n_1''=1$. Moreover, let $\lambda-r=2\alpha+\delta$
with $\delta=0$ or $1$. For $n_2=0$, $\mathfrak{C}_2^*$ vanishes unless
$m_1q_1^2\equiv m_2q_2^2
(\text{mod}\,p^{\alpha})$.
For $n_2\neq 0$, we have
\bna
\mathfrak{C}^*\ll
\widehat{q_1}\widehat{q_2}(\widehat{q_1},\widehat{q_2},n_2)
p^{5(\lambda-r)/2+\min\{k,\alpha\}+3\delta/2}.
\ena
\end{lemma}

\begin{lemma}
Assume $\lambda\leq 2\kappa/3$. Let $p^k\| n_2$ with $k\geq 0$.

(1) We have $\mathfrak{B}^*$ vanishes unless $n_1''=1$ and
\bna
\mathfrak{B}^*\ll \widehat{q_1}\widehat{q_2}(\widehat{q_1},\widehat{q_2},n_2)
\widehat{\rho_s}^2p^{2\lambda}.
\ena
Moreover, for $n_2=0$, the character sums vanish unless $q_1=q_2$ and
$a_1\equiv a_2 (\bmod p^s)$ in which case
\bna
\mathfrak{B}^*\ll \widehat{q_1}^2(\widehat{q_1},m_1-m_2)\widehat{\rho_s}p^{2\lambda+s}.
\ena

(2) Let $\lambda=2\alpha+\delta$
with $\delta=0$ or $1$. For $n_2=0$, we have $\mathfrak{B}_2^*$ vanishes unless
$q_1^2m_1/p^s\equiv
q_2^2m_2/p^s(\text{{\rm mod }} p^{\alpha})$. For $n_2\neq 0$, we have
\bna
\mathfrak{B}^*\ll
\widehat{q_1}\widehat{q_2}(\widehat{q_1},\widehat{q_2},n_2)
p^{5\lambda/2+4s+\min\{k,\alpha\}+3\delta/2}.
\ena

\end{lemma}

For $r\geq \lambda-1$, by (3.15) and Lemma 3 (1), we have
\bea
  \mathscr{H}_1^{\pm}(N,L_1,r)&\ll&N^{\varepsilon}\frac{1}{Q^3}\sum_{n_1'}\sum_{n_1''|p^{\lambda-r}}
  \sum_{1\leq q_1\leq Q\atop n_1'|q_1} \frac{1}{q_1}
 \sum_{1\leq q_2\leq Q\atop n_1'|q_2} \frac{1}{q_2}
        \sum_{1\leq|m_1|\leq N^{\varepsilon}Qp^{\kappa}/N \atop (m_1,q_1)=1}\nonumber\\&&
         \sum_{1\leq|m_2|\leq N^{\varepsilon}Qp^{\kappa}/N \atop (m_2,q_2)=1}
 \frac{1}{\widehat{q_1}\widehat{q_2}\widehat{p_r}}
\sum_{1\leq |n_2|\leq N^{\varepsilon}Q^2p^{\lambda-r}n_1''/L_1}
\widehat{q_1}\widehat{q_2}(\widehat{q_1},\widehat{q_2},n_2)\widehat{p_r}^2p^{2(\lambda-r)}\nonumber\\
&&+N^{\varepsilon}\frac{1}{Q^3}\sum_{n_1'}\sum_{n_1''|p^{\lambda-r}}
  \sum_{1\leq q_1\leq Q,(q_{1},p)=1 \atop n_1'|q_1} \frac{1}{q_1^2}
        \sum_{1\leq|m_1|\leq N^{\varepsilon}Qp^{\kappa}/N\atop (m_1,q_1)=1}
    \sum_{1\leq|m_2|\leq N^{\varepsilon}Qp^{\kappa}/N \atop (m_2,q_1)=1}\nonumber\\&&
    \frac{1}{\widehat{q_1}^2\widehat{p_r}}
\widehat{q_1}^2(\widehat{q_1},m_1-m_2)\widehat{p_r}p^{2(\lambda-r)}\nonumber\\
&\ll&N^{\varepsilon}\left(\frac{Qp^{2\kappa+4\lambda-4r}}{N^2L_1}
+\frac{p^{2\kappa+2\lambda-2r}}{QN^2}
\right).
\eea
For $r\leq \lambda-2$, by Lemma 3, we have
\bea
 \mathscr{H}_1^{\pm}(N,L_1,r)\ll \mathbf{R}_1+\mathbf{R}_2,
\eea
where $\mathbf{R}_1$ is the contribution from $p^{\lambda-r}/n_1''\geq p^2$ and $n_2=0$
\bea
\mathbf{R}_1&=&N^{\varepsilon}\frac{1}{Q^3}\sum_{\delta=0,1}\sum_{n_1'}
  \sum_{1\leq q_1\leq Q \atop n_1'|q_1} \frac{1}{q_1^2}
 \sum_{1\leq|m_1|\leq N^{\varepsilon}Qp^{\kappa}/N \atop (m_1,q_1)=1}
    \sum_{1\leq|m_2|\leq N^{\varepsilon}Qp^{\kappa}/N\atop
    m_2\equiv m_1
\left(\text{mod}\,p^{(\lambda-r-\delta)/2}\right) }\nonumber\\&&
 \frac{1}{\widehat{q_1}^2p^{\lambda-r}}
\widehat{q_1}^2(\widehat{q_1},m_1-m_2)p^{3(\lambda-r)}\nonumber\\
&\ll&N^{\varepsilon}\frac{p^{2(\lambda-r)}}{Q^3}\frac{Qp^{\kappa}}{N}\sum_{\delta=0,1}
\left(1+\frac{Qp^{\kappa}}{N}
p^{-\frac{\lambda-r-\delta}{2}}\right)\nonumber\\
&\ll&N^{\varepsilon}\left(\frac{p^{\kappa+2\lambda-2r}}{Q^2N}+
p^{1/2}\frac{p^{2\kappa+3\lambda/2-3r/2}}{QN^2}\right)
\eea
and $\mathbf{R}_2$ is the remaining piece
\bea
\mathbf{R}_2&=&N^{\varepsilon}\frac{1}{Q^3}\sum_{\delta=0,1}\sum_{n_1'}
  \sum_{1\leq q_1\leq Q\atop n_1'|q_1} \frac{1}{q_1}
 \sum_{1\leq q_2\leq Q\atop n_1'|q_2} \frac{1}{q_2}
    \sum_{1\leq|m_1|\leq N^{\varepsilon}Qp^{\kappa}/N \atop (m_1,q_1)=1}
    \sum_{1\leq|m_2|\leq N^{\varepsilon}Qp^{\kappa}/N\atop (m_2,q_2)=1 }\nonumber\\&&
 \frac{1}{\widehat{q_{1}}\widehat{q_{2}}p^{\lambda-r}}\sum_{0\leq k\ll \log \mathfrak{q}}
\sum_{1\leq|n_2|\leq N^{\varepsilon}Q^2p^{\lambda-r}/L_1 \atop p^k|n_2}
\widehat{q_1}\widehat{q_2}(\widehat{q_1},\widehat{q_2},n_2)
p^{5(\lambda-r)/2+\min\{k,(\lambda-r-\delta)/2\}+3\delta/2}\nonumber\\
&\ll&N^{\varepsilon}\frac{p^{3(\lambda-r)/2}}{Q^3}
\left(\frac{Qp^{\kappa}}{N}\right)^2\sum_{\delta=0,1}
p^{3\delta/2}
\sum_{n_1'}\frac{1}{n_1'}
  \sum_{1\leq q_1\leq Q \atop n_1'|q_1} \frac{1}{q_1}
\left\{\sum_{0\leq k\leq(\lambda-r-\delta_1)/2}\frac{Q^2p^{\lambda-r}}{L_1}\right.\nonumber\\&&\left.
+\sum_{(\lambda-r-\delta_1)/2<k\ll \log \mathfrak{q}}p^{(\lambda-r-\delta)/2}
\frac{Q^2p^{\lambda-r-k}}{L_1}
\right\}\nonumber\\
&\ll&p^{3/2}N^{\varepsilon}\frac{p^{2\kappa+2\lambda-5r/2}}{N^{3/2}L_1}.
\eea
Obviously, the second term in (3.20) is dominated by (3.21),
since $N^{\varepsilon}Q^2p^{\lambda-r}/L_1\geq 1$.
By (3.9) and (3.18-3.21), the contribution from $\mathscr{H}_1^{\pm}(N,L_1,r)$
to $\mathscr{S}(N)$ in (3.6) is at most
\bea
&&N^{\varepsilon}\sum_{r=\lambda-1}^{\lambda}
\sum_{L_1\ll N^{\varepsilon}p^{3\lambda-3r}Q^3/N\atop L_1 \, \mathrm{dyadic}}
\frac{N^{3/2}L_1^{1/2}}{p^{(\kappa+3\lambda-r)/2}}
\left(\frac{Q^{1/2}p^{\kappa+2\lambda-2r}}{NL_1^{1/2}}+\frac{p^{\kappa+\lambda-r}}{Q^{1/2}N}
\right)\nonumber\\
&&+N^{\varepsilon}\sum_{r=0}^{\lambda-2}\sum_{L_1\ll N^{\varepsilon}p^{3\lambda-3r}Q^3/N\atop L_1 \, \mathrm{dyadic}}
\frac{N^{3/2}L_1^{1/2}}{p^{(\kappa+3\lambda-r)/2}}
\left\{\frac{p^{\kappa/2+\lambda-r}}{QN^{1/2}}+
p^{3/4}\frac{p^{\kappa+\lambda-5r/4}}{N^{3/4}L_1^{1/2}}\right\}\nonumber\\
&\ll&N^{1/2+\varepsilon}\sum_{r=\lambda-1}^{\lambda}
\left(N^{1/4}p^{\kappa/2+\lambda/4-3r/2}+p^{\kappa/2+\lambda/2-2r}
\right)\nonumber\\
&+&N^{1/2+\varepsilon}\sum_{r=0}^{\lambda-2}
\left(p^{3\lambda/4-2r}N^{1/4}+p^{3/4}
N^{1/4}p^{\kappa/2-\lambda/2-3r/4}\right)\nonumber\\
&\ll&N^{1/2+\varepsilon}
\left(N^{1/4}p^{3\lambda/4}+p^{3/4}
N^{1/4}p^{\kappa/2-\lambda/2}\right)
\eea
for $\lambda\geq 2$.

Similarly, by Lemma 4, the contribution from $n_2= 0$ to $\mathscr{H}_2^{\pm}(N,L_2,s)$ is at most
\bea
&&N^{\varepsilon} \frac{p^{3s}}{Q^{3}} \sum_{\delta=0,1}\sum_{n_{1}'}
   \sum_{1\leq q_{1}\leq Q/p^{s}\atop n_{1}'|q_{1}} \frac{1}{q_1^2}
    \sum_{1\leq|m_1|\leq N^{\varepsilon}Qp^{\kappa}/N
    \atop p^s|m_1}
\sum_{1\leq|m_2|\leq N^{\varepsilon}Qp^{\kappa}/N,p^s|m_2
    \atop m_2/p^s\equiv m_1/p^s(\text{{\rm mod }}
    p^{(\lambda-\delta)/2})}\nonumber\\&&\qquad\qquad \times
    \frac{1}{\widehat{q_{1}}^2p^{\lambda+s}}
    \widehat{q_{1}}^2\left(\widehat{q_{1}},m_1-m_2\right)p^{3\lambda+2s}\nonumber\\
&=& N^{\varepsilon} \frac{p^{2\lambda+4s}}{Q^{3}}\sum_{\delta=0,1}\sum_{n_1'}
    \mathop{\sum_{1\leq q_{1}\leq Q/p^{s} \atop (q_1,p)=1}}_{n_1'|q_1} \frac{1}{q_1^2}
    \frac{Qp^{\kappa-s}}{N}
    \left(\widehat{q_1}+\frac{Qp^{\kappa-s-(\lambda-\delta)/2}}{N}\right) \nonumber\\
&=& N^{\varepsilon} \left(\frac{p^{\kappa+2\lambda+3s}}{Q^2N}
    +p^{1/2}\frac{p^{2\kappa+3\lambda/2+2s}}{QN^2}\right).
\eea
The contribution from $n_2\neq 0$ to $\mathscr{H}_2^{\pm}(N,L_2,s)$ is bounded by
\bea
&& N^{\varepsilon} \frac{p^{3s}}{Q^{3}} \sum_{\delta=0,1}\sum_{n_{1}'}
    \mathop{\sum_{1\leq q_{1}\leq Q/p^{s} }}_{n_{1}'|q_{1}} \frac{1}{q_{1}}
    \mathop{\sum_{1\leq q_{2}\leq Q/p^{s} }}_{n_{1}'|q_{2}} \frac{1}{q_{2}}
    \sum_{1\leq|m_{1}|\leq N^{\varepsilon}Qp^{\kappa}/N \atop m_{1}\equiv0(\text{{\rm mod }}  p^{s})}
    \sum_{1\leq|m_{2}|\leq N^{\varepsilon}Qp^{\kappa}/N \atop m_{2}\equiv0(\text{{\rm mod }}  p^{s})}
   \nonumber\\&&
\frac{1}{\widehat{q_{1}}\widehat{q_{2}}p^{\lambda+s}}\sum_{0\leq k\ll\log q}
    \sum_{1\leq |n_{2}|\leq N^{\varepsilon}Q^{2}p^{\lambda-s}/L_2\atop p^k|n_2}
    \widehat{q_{1}}\widehat{q_{2}}\left(\widehat{q_{1}},\widehat{q_{1}},n_{2}\right)
    p^{5\lambda/2+4s+\min\{k,(\lambda-\delta)/2\}+3\delta/2}\nonumber\\
&\ll& N^{\varepsilon} \frac{p^{3\lambda/2+6s}}{Q^{3}}\sum_{\delta=0,1} p^{3\delta/2}\sum_{n_{1}'}\frac{1}{n_1'}
    \mathop{\sum_{1\leq q_{1}\leq Q/p^{s} }}_{n_{1}'|q_{1}} \frac{1}{q_{1}}
    \left(\frac{Qp^{\kappa-s}}{N}\right)^2\frac{Q^{2}p^{\lambda-s}}{L_{2}}\nonumber\\
&\ll&  p^{3/2}N^{\varepsilon} \frac{Qp^{2\kappa+5\lambda/2+3s}}{N^{2}L_{2}}.
\eea
Obviously, the second term in (3.23) is bounded by (3.24).
Plugging these estimates into (3.16) and (3.12), we have that the contribution from
$\mathscr{S}_{2}^{\pm}(N,L_{2},s)$ to $\mathscr{S}(N)$ in (3.6) is bounded by
\bea
&&N^{\varepsilon}\sum_{s=1}^{\log Q/\log p}
\sum_{L_2\ll N^{\varepsilon}p^{3\lambda}Q^3/N\atop
L_2 \, \mathrm{dyadic}} \frac{N^{3/2}L_{2}^{1/2}}{p^{(\kappa+3\lambda+3s)/2}}
      \left(\frac{p^{(\kappa+2\lambda+3s)/2}}{QN^{1/2}}+p^{3/4}
      \frac{Q^{1/2}p^{\kappa+5\lambda/4+3s/2}}{NL_2^{1/2}}\right)\nonumber\\
&\ll& N^{1/2+\varepsilon}(N^{1/4}p^{3\lambda/4}
+p^{3/4}N^{1/4}p^{\kappa/2-\lambda/2}).
\eea

\subsection{Conclusion}

By (3.22) and (3.25) we have
\bna
\mathscr{S}(N)\ll N^{1/2+\varepsilon}(N^{1/4}p^{3\lambda/4}
+p^{3/4}N^{1/4}p^{\kappa/2-\lambda/2}).
\ena
Since $N\leq p^{3\kappa/2+\varepsilon}$, Proposition 1 follows.

\section{Proof of Lemma 1}
\setcounter{equation}{0}
\medskip

In this section we apply Poisson summation formula to prove Lemma 1. Precisely,
\bea
\mathscr{A}&=&\sum_{\beta(\mathrm{mod}\,qp^{\kappa})} \chi(\beta) e\left(-\frac{(\overline{a}+bq)\beta}{qp^{\lambda}}\right)
   \sum_{m\equiv\beta(\mathrm{mod}\,qp^{\kappa})} U\left(\frac{m}{N}\right)
   e\left(\frac{m\zeta}{aqp^{\lambda}}\right)\nonumber\\
&:=&\frac{N}{qp^{\kappa}} \sum_{m\in\mathbb{Z}} \mathfrak{A}(m,a,b,q)\mathfrak{I}(m,a,q,\zeta)
\eea
where
\bna
\mathfrak{A}(m,a,b,q)=\sum_{\beta(\mathrm{mod} \,qp^{\kappa})} \chi(\beta)
e\left(\frac{m-(\overline{a}+bq)p^{\kappa-\lambda}}{qp^{\kappa}}\beta\right)
\ena
and
\bea
\mathfrak{I}(m,a,q,\zeta)
=\int_{\mathbb{R}}U(y)e\left(\frac{\zeta Ny}{aqp^{\lambda}}\right)e\left(-\frac{mNy}{qp^{\kappa}}\right)\mathrm{d}y.
\eea

\subsection{Computing the character sum $\mathfrak{A}(m,a,b,q)$}

Write $q=p^{s}q', (q',p)=1$ and $s\geq 0$.
Then
\bna
\mathfrak{A}(m,a,b,q)&=&\sum_{\beta(\mathrm{mod}\, q'p^{s+\kappa})}\chi(\beta)
e\left(\frac{m-(\overline{a}+bq'p^s)p^{\kappa-\lambda}}{q'p^{s+\kappa}}\beta\right)\\
&=&\sum_{\beta_1(\mathrm{mod}\,q')}
   e\left(\frac{m-\overline{a}p^{\kappa-\lambda}}{q'}\overline{p^{s+\kappa}}\beta_1\right)\\&&
\sum_{\beta_2(\mathrm{mod}\,p^{s+\kappa})}\chi(\beta_2)
   e\left(\frac{m-(\overline{a}+bq'p^s)p^{\kappa-\lambda}}{p^{s+\kappa}}\overline{q'}\beta_2\right),
\ena
where the first sum over $\beta_1$ is $q'\mathbf{1}_{m\equiv \overline{a}p^{\kappa-\lambda}(\mathrm{mod}\, q')}$,
and the second sum over $\beta_2$ is
\bna
&&\chi(q')\sum_{\beta_2(\mathrm{mod}\,p^{s+\kappa})}\chi(\beta_2)
   e\left(\frac{m-(\overline{a}+bq'p^s)p^{\kappa-\lambda}}{p^{s+\kappa}}\beta_2\right)\\
&=&\chi(q')\sum_{\alpha_{1}(\mathrm{mod}\,p^{s})}
  e\left(\frac{m-\overline{a}p^{\kappa-\lambda}}{p^{s}}\alpha_1\right)
   \sum_{\alpha_{2}(\mathrm{mod}\,p^{\kappa})}\chi(\alpha_{2})
   e\left(\frac{m-(\overline{a}+bq'p^s)p^{\kappa-\lambda}}{p^{s+\kappa}}\alpha_2\right)\\
&=&p^{s}\chi(q')\mathbf{1}_{m\equiv\overline{a}p^{\kappa-\lambda}(\mathrm{mod}\,p^{s})}
   \overline{\chi}\left(\frac{m-(\overline{a}+bq'p^{s})p^{\kappa-\lambda}}{p^{s}}\right)\tau_{\chi}.
\ena
Here $\tau_{\chi}$ is the Gauss sum. Thus
\bea
\mathfrak{A}(m,a,b,q)=q\mathbf{1}_{m\equiv\overline{a}p^{\kappa-\lambda}(\mathrm{mod}\,q)}\chi(q')
\overline{\chi}\left(\frac{m-(\overline{a}+bq)p^{\kappa-\lambda}}{p^{s}}\right)\tau_{\chi}.
\eea

\subsection{ Bounding the integral $\mathfrak{I}(m,a,q,\zeta)$}
Integration by parts $j$ times, we get
\bna
\mathfrak{I}(m,a,q,\zeta)
\ll\left(\frac{Qp^{\kappa}}{|m|N}\right)^{j}.
\ena
Thus the $m$-sum is essentially supported on $|m|\leq N^{\varepsilon}Qp^{\kappa}/N$.
Then Lemma 1 follows from (4.1) and (4.2).

Furthermore, we also do repeated partial integration by integrating all the exponential factors and
differentiating $U$ only to get
\bna
\mathfrak{I}(m,a,q,\zeta)
\ll\left(\frac{N}{aqp^{\lambda}}\left|\zeta-\frac{ma}{p^{\kappa-\lambda}}\right|\right)^{-j}.
\ena
This restricts the $\zeta$-integral essentially over
$\left|\zeta-ma/p^{\kappa-\lambda}\right|\leq N^{\varepsilon}q/Q$ for any $\varepsilon>0$.

\section{Proof of Lemma 2}
\setcounter{equation}{0}
\medskip
In this section we will apply the $GL_3$ Voronoi formula to transform $\mathscr{B}$, where
\bna
\mathscr{B}=
\sum_{n}A_{\pi}(1,n) e\left(\frac{a^*n}{q^*}\right)\phi(n),
\ena
where $a^*=(\overline{a}+bq)/(\overline{a}+bq,qp^{\lambda})$,
$q^*=qp^{\lambda}/(\overline{a}+bq,qp^{\lambda})$ and $\phi(y)=V\left(y/N\right)e\left(-\zeta y/aqp^{\lambda}\right)$.
Applying the $GL_3$ Voronoi formula (see \cite{GL1}, \cite{MS}), we have
\bea
\mathscr{B}=q^*\sum_{\pm}\sum_{n_{1}|q^*} \sum_{n_{2}=1}^{\infty} \frac{A_{\pi}(n_{2},n_{1})}{n_{1}n_{2}}
  S\left(\overline{a^*},\pm n_{2};\frac{q^*}{n_{1}}\right)
  \Phi_{\phi}^{\pm}\left(\frac{n_{1}^{2}n_{2}}{q^{*3}}\right),
\eea
where
\bna
\Phi_{\phi}^{\pm}\left(y\right)
=\frac{1}{2\pi i}\int\limits_{(\sigma)}y^{-s}
\gamma_{\pm}(s)\widetilde{\phi}(-s)\mathrm{d}s,\qquad \sigma>\max\limits_{1\leq j\leq 3}\{-1-\mathrm{Re}(\mu_j)\},
\ena
where $\mu_j$, $j=1,2,3$, are the Langlands parameters of $\pi$, $\widetilde{\phi}$ is the Mellin transform of $\phi$ and
\bna
\gamma_{\pm}(s)=\frac{1}{2\pi^{3(s+1/2)}}\left(\prod_{j=1}^3
\frac{\Gamma\left((1+s+\mu_j)/2\right)}
{\Gamma\left((-s-\mu_j)/2\right)}\mp i\prod_{j=1}^3
\frac{\Gamma\left((2+s+\mu_j)/2\right)}
{\Gamma\left((-s-\mu_j+1)/2\right)}\right).
\ena

First, we study the integral transform in (5.1). By Stirling's formula, for $\sigma\geq -1/2$,
\bna
\gamma_{\pm}(\sigma+i\tau)\ll_{\pi,\sigma}(1+|\tau|)^{3\left(\sigma+1/2\right)}.
\ena
Moreover, for $s=\sigma+i\tau$,
$$
\widetilde{\phi}(-s)=N^{-s} \widetilde{V}\left(\frac{\zeta N}{aqp^{\lambda}},-s\right)\ll
N^{-\sigma}\min \left\{1,\left(\frac{Q}{q|\tau|}\right)^{j}\right\},
$$
for any $j\geq 0$, where $\widetilde{V}(r,s)=\int_{0}^{\infty} V(y) e(-ry) y^{s-1} \mathrm{d}y$. Thus
\bna
\Phi_{\phi}^{\pm}\left(y\right)\ll \left(\frac{Q}{q}\right)^{5/2}\left(\frac{q^3Ny}{Q^3}\right)^{-\sigma}.
\ena
Thus $ \Phi_{\phi}^{\pm}\left(n_{1}^{2}n_{2}/q^{*3}\right)$ on
the right hand side of (5.1) gives arbitrary power saving in $\mathfrak{q}$ if $n_1^2n_2\geq N^{\varepsilon}q^{*3}Q^3/q^3N$
for any $\varepsilon>0$. For small values of $n_1^2n_2$, we move the integration line to
$\sigma=-1/2$ to get
\bna
\Phi_{\phi}^{\pm}\left(\frac{n_{1}^{2}n_{2}}{q^{*3}}\right)
:=\left(\frac{Nn_{1}^{2}n_{2}}{q^{*3}}\right)^{1/2}\mathfrak{J}^{\pm}
\left(\frac{n_{1}^{2}n_{2}}{q^{*3}},a,q,\zeta\right),
\ena
where
\bea
\mathfrak{J}^{\pm}
\left(y,a,q,\zeta\right)=\frac{1}{2\pi}\int_{\mathbb{R}}(Ny)^{-i\tau}
\gamma_{\pm}\left(-\frac{1}{2}+i\tau\right) \widetilde{V}\left(\frac{\zeta N}{aqp^{\lambda}},\frac{1}{2}-i\tau\right)\mathrm{d}\tau.
\eea
Therefore,
\bea
\mathscr{B}&=&\frac{N^{1/2}}{q^{*1/2}}\sum_{\pm}\sum_{n_{1}|q^*}
\sum_{n_1^2n_2\leq N^{\varepsilon}q^{*3}Q^3/q^3N}\frac{A_{\pi}(n_{2},n_{1})}{\sqrt{n_2}}
  S\left(\overline{a^*},\pm n_{2};\frac{q^*}{n_{1}}\right)\nonumber\\&&\qquad\qquad\qquad\qquad
  \qquad\qquad\qquad
\times\mathfrak{J}^{\pm}\left(\frac{n_{1}^{2}n_{2}}{q^{*3}},a,q,\zeta\right)+O(\mathfrak{q}^{-2018}).
\eea

Furthermore, by the second derivative test for exponential integrals,
\bna
\widetilde{V}\left(\frac{\zeta N}{aqp^{\lambda}},\frac{1}{2}-i\tau\right)\ll (1+|\tau|)^{-1/2}.
\ena
It follows that
\bea
\mathfrak{J}^{\pm}
\left(y,a,q,\zeta\right)\ll N^{\varepsilon}\sqrt{\frac{Q}{q}}.
\eea
Lemma 2 follows from (5.3) and (5.4).

\medskip

\section{Character sums}
\setcounter{equation}{0}
\medskip

In this section we estimate the character sums in (3.14) and (3.17)
\bna
\mathfrak{C}^*=\sum_{\beta(\text{{\rm mod }}\widehat{q_1}\widehat{q_2}\widehat{p_r})}   \mathfrak{C}_r\left(m_{1},n_{1}',n_{1}'',\beta,a_{1},q_{1}\right)
      \overline{\mathfrak{C}_r\left(m_{2},n_1',n_1'',\beta,a_2,q_{2}\right)}
      e\left(\frac{n_2\beta}{\widehat{q_1}\widehat{q_2}\widehat{p_r}}\right)
\ena
and
\bna
\mathfrak{B}^*=\sum_{\beta(\text{{\rm mod }}\widehat{q_1}\widehat{q_2}\widehat{\rho_s})}   \mathfrak{B}_s\left(m_{1},n_{1}',n_{1}'',\beta,a_{1},q_{1}\right)
      \overline{\mathfrak{B}_s\left(m_{2},n_1',n_1'',\beta,a_2,q_{2}\right)}
      e\left(\frac{n_2\beta}{\widehat{q_1}\widehat{q_2}\widehat{\rho_s}}\right),
\ena
where $\widehat{q}=q/n_1'$, $\widehat{p_r}=p^{\lambda-r}/n_{1}''$,
$\widehat{\rho_s}=p^{\lambda+s}/n_1''$, $n_2\in \mathbb{Z}$,
$\mathfrak{C}_r(m,n_1',n_1'',\beta,a,q)$ and $\mathfrak{B}_s(m,n_1',n_1'',\beta,a,q)$ are defined
in (3.11) and (3.13), respectively.
Write $\beta=\widehat{q_{1}}\widehat{q_{2}}\overline{\widehat{q_{1}}\widehat{q_{2}}}b_{1}+\widehat{p_r}
\overline{\widehat{p_r}}b_{2}$, where $b_1\left(\text{{\rm mod }}\widehat{p_r}\right)$ and
$b_2\left(\text{{\rm mod }}\widehat{q_1}\widehat{q_2}\right)$.
We obtain
\bna
\mathfrak{C}^*=\mathfrak{C}_1^*\mathfrak{C}_2^*,
\ena
where
\bna
\mathfrak{C}_1^*=\sum_{b\left(\text{{\rm mod }} \widehat{q_1}\widehat{q_2}\right)}
   S\left(a_1\overline{\varpi_{q_1}^r\widehat{p_r}},b\overline{\widehat{p_r}};\widehat{q_1}\right)
   S\left(a_2\overline{\varpi_{q_2}^r\widehat{p_r}},b\overline{\widehat{p_r}};\widehat{q_2}\right)
   e\left(\frac{n_2\overline{\widehat{p_r}}b}{\widehat{q_{1}}\widehat{q_2}}\right)
\ena
and
\bna
\mathfrak{C}_2^*&=&\sum_{b\left(\text{{\rm mod }}\widehat{p_r}\right)} \sum_{c_1(\text{{\rm mod }} p^{\lambda-r})}
   \overline{\chi}\left(m_1-c_1p^{\kappa-\lambda+r}\right)
   S\left(\overline{c_1}\overline{\widehat{q_{1}}},
 b\overline{\widehat{q_1}};\widehat{p_r}\right)\\
&& \sum_{c_2(\text{{\rm mod }} p^{\lambda-r})}
   \chi\left(m_2-c_2p^{\kappa-\lambda+r}\right)
   S\left(\overline{c_2}\overline{\widehat{q_{2}}},b
   \overline{\widehat{q_2}};\widehat{p_r}\right)
   e\left(\frac{\overline{\widehat{q_{1}}\widehat{q_2}}bn_2}{\widehat{p_r}}\right).
\ena
Similarly,
\bna
\mathfrak{B}^*=\mathfrak{B}_1^*\mathfrak{B}_2^*,
\ena
where
\bna
\mathfrak{B}_1^*=\sum_{b\left(\text{{\rm mod }} \widehat{q_1}\widehat{q_2}\right)}
   S\left(a_1\overline{\widehat{\rho_s}},b\overline{\widehat{\rho_s}};\widehat{q_1}\right)
   S\left(a_2\overline{\widehat{\rho_s}},b\overline{\widehat{\rho_s}};\widehat{q_2}\right)
   e\left(\frac{n_2\overline{\widehat{\rho_s}}b}{\widehat{q_{1}}\widehat{q_2}}\right)
\ena
and
\bna
\mathfrak{B}_2^*&=&\sum_{b\left(\text{{\rm mod }}\widehat{\rho_s}\right)} \sum_{c_1(\text{{\rm mod }} p^{\lambda})}
   \overline{\chi}\left(\frac{m_1-(\overline{a_1}+c_1p^s)p^{\kappa-\lambda}}{p^s}\right)
   S\left(\overline{\overline{a_1}+c_1p^s}\overline{\widehat{q_{1}}},
 b\overline{\widehat{q_1}};\widehat{\rho_s}\right)\\
&& \sum_{c_2(\text{{\rm mod }} p^{\lambda})}
   \chi\left(\frac{m_{2}-(\overline{a_{2}}+c_2p^s)p^{\kappa-\lambda}}{p^s}\right)
   S\left(\overline{\overline{a_{2}}+c_2p^s}\overline{\widehat{q_{2}}},b
   \overline{\widehat{q_2}};\widehat{\rho_s}\right)
   e\left(\frac{\overline{\widehat{q_{1}}\widehat{q_2}}bn_2}{\widehat{\rho_s}}\right).
\ena

We quote the following estimates for $\mathfrak{C}_1^*$ and $\mathfrak{B}_1^*$
which were proved in \cite{Munshi31}
by induction.

\begin{lemma}
We have
\bna
\mathfrak{C}_1^*, \mathfrak{B}_1^* \ll \widehat{q_1}\widehat{q_2}(\widehat{q_1},\widehat{q_2},n_2).
\ena
Moreover, for $n_2=0$, the character sums vanish unless $q_1=q_2$ in which case
\bna
\mathfrak{C}_1^*,  \mathfrak{B}_1^*\ll \widehat{q_1}^2(\widehat{q_1},m_1-m_2).
\ena

\end{lemma}

For $\mathfrak{C}_2^*$ and $\mathfrak{B}_2^*$, we will prove the following results.

\begin{lemma}
Assume $\lambda\leq 2\kappa/3$. Let $p^k\| n_2$ with $k\geq 0$.

(1) We have
\bea
\mathfrak{C}_2^*\ll \widehat{p_r}^2p^{2(\lambda-r)}.
\eea
Moreover, for $n_2=0$, we have
\bea
\mathfrak{C}_2^*\ll \widehat{p_r}p^{2(\lambda-r)}.
\eea

(2) For $n_1''=p^{\lambda-r}$ or $n_1''=p^{\lambda-r-1}$ with $\lambda-r\geq 2$, we have
\bna
\mathfrak{C}_2^*=0.
\ena

(3) For $p^{\lambda-r}/n_1''\geq p^2$, we have
$\mathfrak{C}_2^*$ vanishes unless $n_1''=1$. Moreover, let $\lambda-r=2\alpha+\delta$
with $\delta=0$ or $1$. For $n_2=0$, $\mathfrak{C}_2^*$ vanishes unless
$m_1q_1^2\equiv m_2q_2^2
(\text{mod}\,p^{\alpha})$.
For $n_2\neq 0$, we have
\bea
\mathfrak{C}_2^*\ll
p^{5(\lambda-r)/2+\min\{k,\alpha\}+3\delta/2}.
\eea

\end{lemma}

\begin{lemma}
Assume $\lambda\leq 2\kappa/3$. Let $p^k\| n_2$ with $k\geq 0$.

(1) We have $\mathfrak{B}_2^*$ vanishes unless $n_1''=1$ and
\bna
\mathfrak{B}_2^*\ll \widehat{\rho_s}^2p^{2\lambda}.
\ena
Moreover, for $n_2=0$, we have $a_2\equiv\widehat{q_2}^2\overline{\widehat{q_1}^2}a_1
(\text{{\rm mod }} p^s)$
in which case
\bna
\mathfrak{B}_2^*\ll \widehat{\rho_s}p^{2\lambda+s}.
\ena

(2) Let $\lambda=2\alpha+\delta$ with $\delta=0$ or 1.
For $n_2=0$, we have $\mathfrak{B}_2^*$ vanishes unless $q_1^2m_1/p^s\equiv
q_2^2m_2/p^s(\text{{\rm mod }} p^{\alpha})$.
For $n_2\neq0$, we have
\bna
\mathfrak{B}_2^*\ll p^{5\lambda/2+4s+\min\{k,\alpha\}+3\delta/2}.
\ena

\end{lemma}
Now Lemma 3 follows from Lemmas 5 and 6, and Lemma 4 follows from Lemmas 5 and 7.
We only prove Lemma 6 for $\mathfrak{C}_2^*$ in detail,
since the proof of Lemma 7 is very similar.
\begin{proof}
(1) Trivially, (6.1) follows from Weil's bound for Kloosterman sums.
To prove (6.2), we open the Kloosterman sums and sum over $b$ to get
\bea
\mathfrak{C}_2^*
&=&\widehat{p_r} \sum_{c_1(\text{{\rm mod }} p^{\lambda-r})}
   \overline{\chi}\left(m_1-c_1p^{\kappa-\lambda+r}\right)
   \sum_{c_2(\text{{\rm mod }} p^{\lambda-r})}
   \chi\left(m_2-c_2p^{\kappa-\lambda+r}\right)\nonumber\\
&& \sideset{}{^*}\sum_{d\left(\text{{\rm mod }}\widehat{p_r}\right)}
   e\left(\frac{\overline{\widehat{q_2}c_2}
   -\widehat{q_2}
   \overline{\widehat{q_1}\left(\widehat{q_1}+ n_2d\right)
   c_1}}{\widehat{p_r}}d\right).
\eea
For $n_2=0$, we denote $m_0=\overline{\widehat{q_2}c_2}
   -\widehat{q_2}
   \overline{\widehat{q_1}^2c_1}$. Since the Ramanujan sum
   $$S(m,0;c)=\sideset{}{^*}\sum\limits_{\alpha(\text{{\rm mod }} c)}\exp(2\pi im\alpha/c)=
   \mu\left(c/(m,c)\right)\varphi(c)/\varphi\left(c/(m,c)\right),$$
   where $\mu$ is the M\"{o}bius function and $\varphi$ is the Euler function,
   the last sum over $d$ for $n_2=0$ is
\bna
\mu\left(\frac{\widehat{p_r}}{(m_0,\widehat{p_r})}\right)
\frac{\varphi(\widehat{p_r})}{\varphi\left(\widehat{p_r}/(m_0,\widehat{p_r})\right)}
=\left\{\begin{array}{ll}
\widehat{p_r}(1-p^{-1}),&\mbox{if}\,(m_0,\widehat{p_r})=\widehat{p_r}\\
-\widehat{p_r}/p,&\mbox{if}\,(m_0,\widehat{p_r})=\widehat{p_r}/p\\
0,& \mbox{otherwise}.
\end{array}\right.
\ena
Thus (6.2) follows.

(2) Let $\lambda-r=2\alpha+\delta$ with $\delta=0$ or 1, $\alpha\geq 1$ is a positive integer.
Write $c_1=b_1p^{\alpha+\delta}+b_2$, $b_1(\text{{\rm mod }}p^{\alpha})$,
$b_2(\text{{\rm mod }}p^{\alpha+\delta})$ and
$c_2=h_1p^{\alpha+\delta}+h_2$, $h_1(\text{{\rm mod }}p^{\alpha})$,
$h_2(\text{{\rm mod }}p^{\alpha+\delta})$. If $n_1''=p^{\lambda-r-1}$, we have $\widehat{p_r}=p$ and
\bna
\mathfrak{C}_2^*
&=&p\sum_{b_2(\text{{\rm mod }} p^{\alpha+\delta})}
\sum_{h_2(\text{{\rm mod }} p^{\alpha+\delta})}
\overline{\chi}\left(m_1-b_2p^{\kappa-2\alpha-\delta}\right)
\chi\left(m_2-h_2p^{\kappa-2\alpha-\delta}\right)
\nonumber\\
   &&\sideset{}{^*}\sum_{d\left(\text{{\rm mod }}p\right)}
   e\left(\frac{\overline{\widehat{q_2}h_2}
   -\widehat{q_2}
   \overline{\widehat{q_1}\left(\widehat{q_1}+ n_2d\right)
   b_2}}{p}d\right)
\sum_{b_1(\text{{\rm mod }} p^{\alpha})}
\chi\left(1+\overline{m_{1}-b_2
   p^{\kappa-2\alpha-\delta}} p^{\kappa-\alpha}b_1\right)\nonumber\\&&
\sum_{h_1(\text{{\rm mod }} p^{\alpha})}
 \chi\left(1-\overline{m_{2}-h_2p^{\kappa-2\alpha}} p^{\kappa-\alpha}h_1\right).
\ena
Recall $\chi$ is a primitive character of modulus $p^{\kappa}$ and $\kappa>\lambda\geq 2\alpha$.
Thus $\chi(1+zp^{\kappa-\alpha})$
is an additive character to modulus $p^{\alpha}$, so there exists an integer
$\eta$ (uniquely determined modulo $p^{\alpha}$), $(\eta,p)=1$, such that
$\chi(1+zp^{\kappa-\alpha})=\exp(2\pi i \eta z/p^{\alpha})$. Therefore,
$\mathfrak{C}_2^*=0$. For $n_1''=p^{\lambda-r}$, the proof is similar and easier.

(3) Write
$p^{\lambda-r}=p^{2\alpha+\delta_1}$, $\widehat{p_r}=p^{\lambda-r}/n_1''=p^{2\beta+\delta_2}$,
$\delta_1=0$ or 1, $\delta_2=0$ or 1,
$\alpha\geq 1$
and $\beta\geq 1$.
Write $c_1=b_1p^{\alpha+\delta_1}+b_2$, $b_1(\text{{\rm mod }}p^{\alpha})$,
$b_2(\text{{\rm mod }}p^{\alpha+\delta_1})$,
$c_2=h_1p^{\alpha+\delta_1}+h_2$, $h_1(\text{{\rm mod }}p^{\alpha})$,
$h_2(\text{{\rm mod }}p^{\alpha+\delta_1})$,
and $d=d_1p^{\beta+\delta_2}+d_2$,
$d_1(\text{{\rm mod }}p^{\beta})$, $d_2(\text{{\rm mod }}p^{\beta+\delta_2})$. Then by (6.4), we have
\bna
\mathfrak{C}_2^*
&=&p^{2\beta+\delta_2}\sum_{b_2(\text{{\rm mod }} p^{\alpha+\delta_1})}\sum_{h_2(\text{{\rm mod }} p^{\alpha+\delta_1})}\,
 \sideset{}{^*}\sum_{d_2\left(\text{{\rm mod }}p^{\beta+\delta_2}\right)}
 \sum_{b_1(\text{{\rm mod }} p^{\alpha})}\sum_{h_1(\text{{\rm mod }} p^{\alpha})}\,
 \sum_{d_1\left(\text{{\rm mod }}p^{\beta}\right)}
   \\&&\overline{\chi}\left(m_1-(b_2+b_1p^{\alpha+\delta_1})p^{\kappa-2\alpha-\delta_1}\right)
   \chi\left(m_2-(h_2+h_1p^{\alpha+\delta_1})p^{\kappa-2\alpha-\delta_1}\right)\\
&&e\left(\frac{
   \overline{\widehat{q_2}(h_2+h_1p^{\alpha+\delta_1})}-\widehat{q_2}
   \overline{\widehat{q_1}\left(\widehat{q_1}+ n_2d_2+n_2d_1p^{\beta+\delta_2}\right)
   (b_2+b_1p^{\alpha+\delta_1})}}{p^{2\beta+\delta_2}}(d_2+d_1p^{\beta+\delta_2})\right).
\ena
Note that $\kappa>\lambda\geq 2\alpha+\delta_1\geq 2\beta+\delta_2$ and
$\overline{a+bp^{\alpha}}\equiv \overline{a}(1-\overline{a}bp^{\alpha})(\text{{\rm mod }}p^{2\alpha})$. Thus
\bea
\mathfrak{C}_2^*
=p^{2\alpha+3\beta+\delta_2} \mathop{\sum_{b_2(\text{{\rm mod }} p^{\alpha+\delta_1})}
\sum_{h_2(\text{{\rm mod }} p^{\alpha+\delta_1})}\,
 \sideset{}{^*}\sum_{d_2\left(\text{{\rm mod }}p^{\beta+\delta_2}\right)}}_{
 \overline{(\widehat{q_1}+n_2d_2)^2b_2}
    \widehat{q_2}^2n_2d_2- \overline{
   (\widehat{q_1}+n_2d_2)b_2}\widehat{q_2}^2
   +\widehat{q_1}\overline{h_2}\equiv 0(\text{{\rm mod }} p^{\beta}) }f(b_2,h_2,d_2)\mathcal {C}_1\mathcal {C}_2,
\eea
where
\bna
\mathcal {C}_1&=&\frac{1}{p^{\alpha}}\sum_{b_1(\text{{\rm mod }} p^{\alpha})}
   \chi\left(1+\overline{m_{1}-b_2
   p^{\kappa-2\alpha-\delta_1}}p^{\kappa-\alpha}b_{1}\right)
   e\left(\frac{\widehat{q_2}\overline{\widehat{q_1}(\widehat{q_1}+n_2d_2)b_2^2}
                d_2n_1''}{p^{\alpha}}b_1\right),\\
\mathcal {C}_2&=&\frac{1}{p^{\alpha}}\sum_{h_1( \text{{\rm mod }}p^{\alpha})}
   \chi\left(1-\overline{m_{2}-h_2p^{\kappa-2\alpha-\delta_1}}p^{\kappa-\alpha}h_1\right)
   e\left(\frac{-\overline{\widehat{q_2}h_2^2} d_2 n_1''}
   {p^{\alpha}}h_1\right)
\ena
and
\bna
f(b_2,h_2,d_2)=\overline{\chi}\left(m_1-b_2p^{\kappa-2\alpha-\delta_1}\right)
   \chi\left(m_2-h_2p^{\kappa-2\alpha-\delta_1}\right)
  e\left(\frac{\overline{\widehat{q_2}h_2}
  -\widehat{q_2}\overline{\widehat{q_1}(\widehat{q_1}+n_2d_2)b_2
  }}{p^{2\beta+\delta_2}}d_2\right).
\ena
Since $\chi(1+zp^{\kappa-\alpha})=\exp(2\pi i \eta z/p^{\alpha})$ with $(\eta,p)=1$, we have
\bna
\mathcal {C}_1&=&\frac{1}{p^{\alpha}}\sum_{b_1(\text{{\rm mod }}  p^{\alpha})}
  e\left(\frac{\overline{m_{1}-b_2
   p^{\kappa-2\alpha-\delta_1}} \eta}{p^{\alpha}}b_1\right)
   e\left(\frac{\widehat{q_2}\overline{\widehat{q_1}(\widehat{q_1}+n_2d_2)b_2^2}
                d_2n_1''}{p^{\alpha}}b_1\right)\\
&=&\mathbf{1}_{\overline{m_{1}-b_2
   p^{\kappa-2\alpha-\delta_1}} \eta+\widehat{q_2}\overline{\widehat{q_1}(\widehat{q_1}+n_2d_2)b_2^2}
                d_2n_1''\equiv 0 (\text{{\rm mod }}p^{\alpha})}.
\ena
Thus $\mathcal {C}_1$ vanishes unless $n_1''=1$ which in turn implies that $\alpha=\beta$ and $\delta_1=\delta_2$.
Moreover, by taking $\lambda\leq 2\kappa/3$, we have $\kappa\geq 3\alpha+2\delta_1$. Hence
$\mathcal {C}_1$ vanishes unless $\overline{m_1} \eta+\widehat{q_2}
\overline{\widehat{q_1}(\widehat{q_1}+n_2d_2)b_2^2}
                d_2\equiv 0 (\text{mod}\,p^{\alpha})$.
Similarly,
\bna
\mathcal {C}_2=\mathbf{1}_{\overline{m_2} \eta+\overline{\widehat{q_2}h_2^2} d_2
   \equiv 0 (\text{{\rm mod }}p^{\alpha})}.
\ena
Plugging these into (6.5) we obtain
\bea
\mathfrak{C}_2^*
&=&p^{5\alpha+\delta_1} \mathop{\sum_{b_2(\text{{\rm mod }} p^{\alpha+\delta_1})}
\sum_{h_2(\text{{\rm mod }} p^{\alpha+\delta_1})}\,
 \sideset{}{^*}\sum_{d_2\left(\text{{\rm mod }}p^{\alpha+\delta_1}\right)}}_{\substack{
 \overline{(\widehat{q_1}+n_2d_2)^2b_2}
    \widehat{q_2}^2n_2d_2- \overline{
   (\widehat{q_1}+n_2d_2)b_2}\widehat{q_2}^2
   +\widehat{q_1}\overline{h_2}\equiv 0(\text{{\rm mod }} p^{\alpha})\\
   \overline{m_1} \eta+\widehat{q_2}\overline{\widehat{q_1}(\widehat{q_1}+n_2d_2)b_2^2}
                d_2\equiv 0 (\text{mod}\,p^{\alpha})\\
   \overline{m_2} \eta+\overline{\widehat{q_2}h_2^2} d_2
   \equiv 0 (\text{mod}\,p^{\alpha})
   }}f(b_2,h_2,d_2).
\eea
To count the numbers of $b_2,h_2$ and $d_2$, we solve the three congruence
equations in (6.6).

(i) If $n_2=0$ or $n_2=p^kn_2'$ with $(n_2',p)=1$ and $p^k\geq p^{\alpha}$, we have
\bna
\left\{\begin{array}{l}
h_2\equiv
\overline{\widehat{q_2}^2}\widehat{q_1}^2b_2
  (\text{{\rm mod }} p^{\alpha})\\
   d_2\equiv -\overline{m_1} \eta\widehat{q_1}^2\overline{\widehat{q_2}}
                b_2^2 (\text{mod}\,p^{\alpha})\\
   d_2\equiv -\overline{m_2} \eta\widehat{q_2}h_2^2
   (\text{mod}\,p^{\alpha})
\end{array}
\right.
\ena
By the last two equations, one sees that $\mathfrak{C}_2^*$ vanishes unless
$m_1\widehat{q_1}^2\equiv m_2\widehat{q_2}^2
(\text{mod}\,p^{\alpha})$. Moreover, for fixed $b_2$, $h_2$ and $d_2$ are uniquely determined
modulo $p^{\alpha}$. Therefore,
\bea
\mathfrak{C}_2^*\ll p^{6\alpha+4\delta_1}
\ll p^{3(\lambda-r)+\delta_1}.
\eea

(ii) If $n_2\neq 0$, we let $n_2=p^kn_2'$ with $(n_2',p)=1$ and $p^k< p^{\alpha}$, and let $\gamma=\overline{\widehat{q_1}+n_2d_2}$. Then
$d_2\equiv\overline{n_2'}
(\overline{\gamma}-\widehat{q_1})/p^k\left(\text{{\rm mod }}  p^{\alpha-k}\right)$ and
the three equations give
\bea
\left\{\begin{array}{l}
b_2\equiv \widehat{q_2}^2\gamma^2h_2(\text{{\rm mod }} p^{\alpha}),\\
\gamma\equiv\overline{\widehat{q_1}}\left(1+\overline{m_1}\eta\widehat{q_1}
\overline{\widehat{q_2}}
n_2b_2^2\right)(\text{{\rm mod }} p^{\alpha}),\\
\overline{\gamma}\equiv
\widehat{q_1}\left(1-\overline{m_2}\eta\overline{\widehat{q_1}}
\widehat{q_2}
n_2h_2^2\right)(\text{{\rm mod }} p^{\alpha}).
\end{array}\right.
\eea
Plugging the second equation into the first equation in (6.8) we get
\bea
b_2\equiv \widehat{q_2}^2\overline{\widehat{q_1}}^2
\left(1+\overline{m_1}\eta\widehat{q_1}\overline{\widehat{q_2}}
n_2b_2^2\right)^2h_2(\text{{\rm mod }}p^{\alpha}).
\eea
By (6.9) and the last two equations in (6.8) we get
\bna
&& \left(\overline{m_{1}}\eta\widehat{q_1}\overline{\widehat{q_2}}\right)^5u^5
   +4\left(\overline{m_1}\eta\widehat{q_1}\overline{\widehat{q_2}}\right)^4u^4
   +6\left(\overline{m_1}\eta\widehat{q_{1}}\overline{\widehat{q_{2}}}\right)^3u^3
   +4\left(\overline{m_1}\eta\widehat{q_1}\overline{\widehat{q_{2}}}\right)^2u^2\\
&& \qquad-\overline{m_1m_2}\eta^2\widehat{q_1}^4\overline{\widehat{q_2}}^4u^2
   +\overline{m_1}\eta\widehat{q_1}\overline{\widehat{q_2}}u
   -\overline{m_2}\eta\widehat{q_1}^3\overline{\widehat{q_2}}^3u
   \equiv0(\text{{\rm mod }} p^{\alpha}),
\ena
where $u=n_2b_2^2$. Thus there are at most 5 roots modulo $p^{\alpha}$ for $u$.
Therefore, there are at most 10 roots modulo $p^{\alpha-k}$ for $b_2$.
For fixed $u$,
$\gamma$ is uniquely determined modulo $p^{\alpha}$ and for fixed $\gamma$ and $b_2$,
$h_2$ is uniquely determined modulo $p^{\alpha}$ by the first equation in (6.8).
Then by the last congruence equation in
(6.6), $d_2$ is uniquely determined modulo $p^{\alpha}$.
Therefore,
\bna
\mathfrak{C}_2^*\ll p^{5\alpha+k+4\delta_1}\ll p^{5(\lambda-r)/2+k+3\delta_1/2}.
\ena
By (6.7) and (6.10), the bound in (6.3) follows.

\end{proof}

\bigskip

\noindent
{\sc Acknowledgements.}
The first author is supported by IRT16R43,
Young Scholars Program of Shandong University, Weihai (Grant No. 2015WHWLJH04) and
Natural Science Foundation of Shandong Province (Grant No. ZR2016AQ15).

\bigskip

\bigskip

{\sc \small School of Mathematics and Statistics, Shandong University, Weihai,
Weihai, Shandong 264209, China}

{\footnotesize {\it E-mail address:} qfsun@sdu.edu.cn}

{\footnotesize {\it E-mail address:} zhaoruisdu@gmail.com}

\end{document}